\documentclass[a4paper,12pt]{article}

\usepackage[english]{babel}
\usepackage[T2A]{fontenc}
\usepackage[utf8]{inputenc}

\usepackage[left=2cm,right=2cm,top=2cm,bottom=2cm]{geometry}
\usepackage{amsthm, amsfonts, amssymb, amsthm, mathtools, dsfont}
\usepackage{comment}

\DeclareMathOperator{\pr}{Pr}
\DeclareMathOperator{\fl}{Fl}
\DeclareMathOperator{\T}{T}
\DeclareMathOperator{\Cy}{C}

\theoremstyle{plain}
\newtheorem{theorem}{Theorem}
\newtheorem{lemma}{Lemma}
\newtheorem*{corollary}{Corollary}

\title{Equivalence of flat-virtual diagrams}
\author{D.A.Popova}

\begin{document}
	\maketitle
	\begin{abstract}
		Maps from links in thickened surfaces to flat-virtual links help to construct invariants of links using invariants of flat-virtual links. This work is dedicated to investigation of equivalence and invariants of flat-virtual diagrams obtained as images of maps from links in thickened surfaces.
	\end{abstract}
	\setlength{\parskip}{0.2cm}
	
	Mathematical Subject Classification: $57M25$
	
	{\em Key words and phrases.} links in thickened surfaces, flat-virtual links, Jones polynomial, torus links

	\section{Introduction}
	
	In paper [1] V.O.Manturov and I.M.Nikonov for a given natural $d$ (lattice $l$, respectively) constructed maps $\phi_d$ ($\phi'_l$, respectively) from links in the full torus (thickened torus, respectively) to flat-virtual links and presented Alexander-like polynomial $\Delta$ and Jones polynomial $J_f$ for flat-virtual diagrams. The Jones polynomial is a picture-valued invariant which means that its values contain equivalence classes of diagrams. In this work we present direct algorithms of constructing the compositions $J_f \circ \phi_d$ and $J_f \circ \phi'_l$, whose values are easily classified.
	
	A \textit{flat-virtual diagram} is a four-valent graph on the plane where each vertex is of one of the following three types: classical (in this case one pair of edges is marked as an overcrossing strand), flat, virtual. In figures we mark flat and virtual crossings as shown in Fig.1. \textit{Flat-virtual link} is an equivalence class of flat-virtual diagrams modulo the following moves:
	
	1. Classical Reidemeister moves;
	
	2. Second and third flat Reidemeister moves which refer to flat crossings only;
	
	3. Third mixed flat-classical Reidemeister move (a strand containing two consecutive flat crossings passes through a classical crossing, the over-under structure at the classical crossing is preserved);
	
	4. Virtual detour moves (a strand containing only virtual crossings can be removed and replaced with a strand with the same endpoints where all new crossing are marked virtual.
	
	\begin{center}
		\includegraphics[width = 1.0\textwidth]{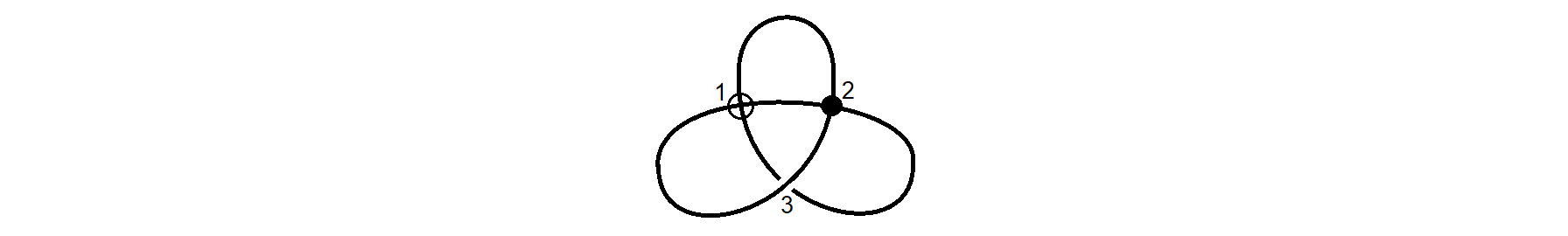}
		Figure 1: A flat-virtual diagram: crossing 1 is virtual, crossing 2 is flat, crossing 3 is classical
	\end{center}
	
	\subsection{The flat-virtual Jones polynomial}
	Let $K$ be flat-virtual diagram. Recall the definition of flat-virtual bracket and flat-virtual Jones polynomial following [1]. Every classical crossing can be smoothed in two ways: A-smoothing (\includegraphics[width=20pt]{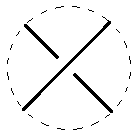} $\rightarrow$ \includegraphics[width=20pt]{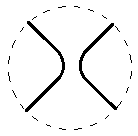}) and B-smoothing (\includegraphics[width=20pt]{classcross} $\rightarrow$ \includegraphics[width=20pt]{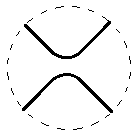}). A {\em state} of a diagram is a set of smoothings of all classical crossings. For a state $s$ $K_s$ is defined as a flat-virtual diagram that we get after putting its crossings in the state $s$. Then flat-virtual bracket $\langle K\rangle_{flat}$ takes value in $Z[a, a^{-1}] \mathcal{G}_f$, where $\mathcal{G}_f$ is a set of equivalence classes of flat-virtual diagrams, and is defined as $\langle K\rangle_{flat} = \displaystyle \sum_s {a^{\alpha(s)-\beta(s)}K_s}$, where $s$ runs over all states of $K$, $\alpha(s)$ ($\beta(s)$ respectively) equals the quantity of A-smoothings (B-smoothings) of the state s. Finally we can define flat virtual Jones polynomial $J_f$ for an oriented flat-virtual link $\overline{K}$ as follows: $J_f(\overline{K}) = (-a)^{-3w(\overline{K})} \langle K\rangle_{flat}$. In [1] it is proven that $J_f$ is an invariant of flat-virtual links.
	
	\subsection{Maps $\phi_d$ and $\phi'_l$}
	Given a cylinder $C = S^1 \times [0, 1]$ with angular coordinate $\alpha \in [0, 2\pi = 0)$ and vertical coordinate
	$z \in [0, 1]$. Fix a number $d \in \mathbb{N}$ and say that two points $x$, $y$ on $C$ are {\em equivalent} if they have the same vertical coordinate and their angular coordinate differs by $\frac{2 \pi k}{d}$ for some integer $k$. Likewise, for the torus $T^2$ with two angular coordinates $(\alpha, \beta) \in [0, 2\pi = 0)$ we fix a lattice $l$ as a discrete subgroup od $T^2$. Say, choosing $d_1, d_2$ our group contains points with angular coordinates $(\frac{2 \pi k_1}{d_1}, \frac{2 \pi k_2}{d_2})$. We say that two points $A$, $B$ on the torus are {\em equivalent} if their difference $A - B$ belongs to the subgroup $l$.
	
	Now having a diagram $K$ on the cylinder (torus respectively) with a fixed number $d$ (lattice $l$, respectively) we define a flat-virtual diagram $\phi_d(K)$ ($\phi'_l$) (up to detour moves) as follows. First, we assume without loss of generality that $K$ is generic with respect to the subgroup, i.e., there is no three pairwise equivalent points, none of the points from the pair of equivalent ones is a crossing and the tangent vectors of the diagram at the equivalent points lying on the edges of the diagram are transverse. In both cases (for $\phi$ and for $\phi'$) for each pair of equivalent points on the edges (say, ($e_j$ , $e'_j$)) we create a new flat crossing in the following manner. Let $e'_j=e_j + g$, where $g$ is an element of the corresponding group ($Z_d$ or $l$). We choose one of these points (never mind which one, say, $e_j$), remove a small neighborhood $U(e_j)$ with endpoints $A$, $B$, shift it by g and connect $A$ to $A + g$ and $B$ to $B + g$ by any curves on the plane not passing through any crossings (all crossings on the newborn curves $[A, A + g]$, $[B, B + g]$ are declared to be virtual).
	
	For a flat-virtual diagram $K'$ with $2m$ cut ends (as shown in Fig.3.2 for $m = 3$) in the rectangle define the flat-virtual diagram $\Cy(K')$, so that the ends $1, ..., m$ go from the right side from top to bottom, the ends $m+1, ..., 2m$ go from the left side from top to bottom: connect $i$ with $m + i$ for $i=1, ..., m$ outside the rectangle without crossings.
	Analogously, for a flat-virtual diagram $K'$ with $2(m_1 + m_2)$ cut ends (as shown in Fig.3.3 for $m_1 = 3$, $m_2 = 2$)in the rectangle define the flat-virtual diagram $\T(K')$, so that the ends $1, ..., m_1$ go from the right side from top to bottom, the ends $m_1 + 1, ..., m_1 + m_2$ go from the top side from left to right, the ends $m_1 + m_2 + 1, ..., 2 m_1 + m_2$ go from the left side from top to bottom, the ends $2m_1 + m_2 + 1, ..., 2(m_1 + m_2)$ go from the bottom side from left to right: connect $i$ with $m_1 + m_2 + i$ for $i = 1, ..., m_1$ and $m_1 + j$ with $2m_1 + m_2 + j$ for $j = 1, ..., m_2$ outside the rectangle, putting all crossings virtual there. Fig.2 shows examples of $\Cy(K')$ and $\T(K')$.
	
	\begin{center}
		\includegraphics[width = 1.0\textwidth]{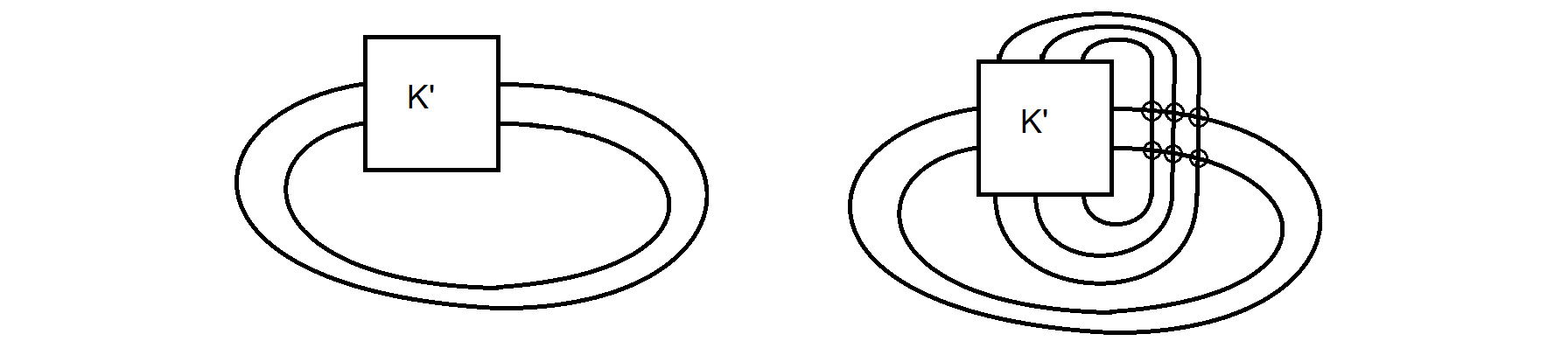}
		Figure 2.1: $\Cy(K')$, $m = 2$, 2.2: $\T(K')$, $m_1 = 3$, $m_2 = 2$
	\end{center}
	
	Split the cylinder into $d$ rectangles with angular width $\frac{2 \pi}{d}$ and height equal to the height of the cylinder. Analogously, split torus into $d_1 d_2$ rectangles with angular width $\frac{2 \pi}{d_1}$ and angular height $\frac{2 \pi}{d_2}$. Notice the following way of constructing $\phi_d(K)$ ($\phi'_l(K)$ respectively) where $K$ is a diagram on cylinder (torus) generic in the following way: all pairs of equivalent points lying on the diagrams edges do not lie on the sides of the rectangles. Let $\pr_d$ ($\pr'_l$) be a following projection map from cylinder (torus) to a rectangle on a plane with width equal to $\frac{2 \pi}{d}$ ($\frac{2 \pi}{d_1}$) and height equal to the height of the cylinder ($\frac{2 \pi}{d_2}$). For a cylinder $c$ and point $X \in c$ with angular coordinate $\alpha + \frac{2 \pi k}{d}$, $\alpha \in [0; \frac{2 \pi}{d})$ and height coordinate $h$  $\pr_d(X) = (\alpha, h)$. For a torus $t$ and point $X \in t$ with angular coordinates ($\alpha_1 + \frac{2 \pi k_1}{d_1}$, $\alpha_2 + \frac{2 \pi k_2}{d_2}$) $\alpha_1 \in [0; \frac{2 \pi}{d_1})$, $\alpha_2 \in [0; \frac{2 \pi}{d_2})$ $\pr_d(X) = (\alpha_1, \alpha_2)$. Then $K' = \pr_d(K)$ ($K' = \pr'_l(K)$) is a flat-virtual diagram on rectangle (we put new crossings flat) with classical and flat crossings and $2m$ ($2(m_1 + m_2)$) cut ends. Then $\phi_d(K) = \Cy(K')$ ($\phi'_l(K) = \T(K')$).
	
	\section{Untangling $\phi_d(K)$ and $\phi'_l(K)$}
	Having maps $\phi_d$ and $\phi'_l$ we shall find out which invariants of flat virtual diagrams distinguish images of the maps well. Let us call a flat-virtual diagram {\em semi-trivial} if it is equivalent to a disjoint union of trivial components and {\em flat eight} components (Fig.3.1). In this section we will see that any image of $\phi_d$ or $\phi'_l$ with smoothed classical crossings is semi-trivial. This means that applying flat virtual Jones polynomial to images of $\phi_d$ and $\phi'_l$ gives sum of polynomials multiplied by disjoint union of flat eight components, which can be easily represented by a two-variable polynomial.
	\begin{center}
    	\includegraphics[width = 1.0\textwidth]{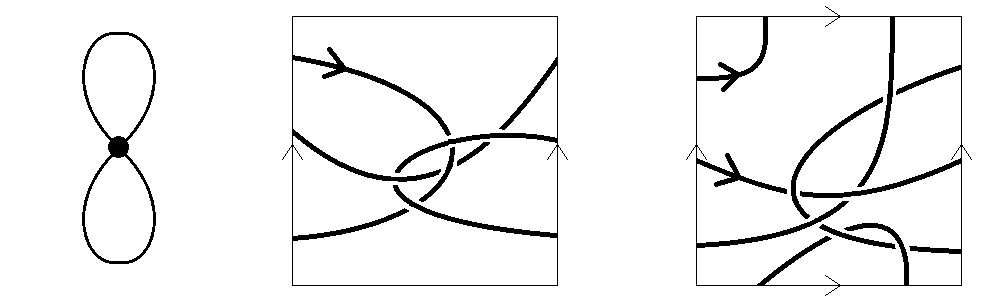}
    	Figure 3.1: Flat eight component, 3.2: Diagram on the cylinder, 3.3 Diagram on the torus
    \end{center}
    The following lemma is useful for untangling the images of mappings $\phi_d$ and $\phi'_l$.
    \begin{center}
    	\includegraphics[width = 1.0\textwidth]{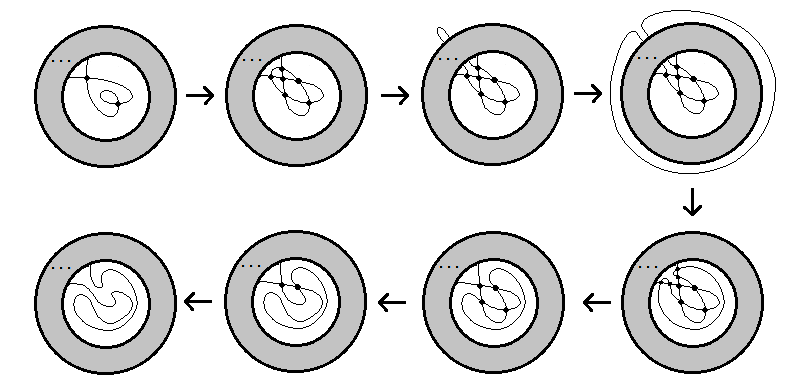}
    	Figure 4: Proof of lemma 1 %стрелки
    \end{center}
    \begin{lemma}
    	Any flat virtual diagram $L$ that is realisable without virtual crossings and has a subdiagram $L'$ as shown in the Fig.4 (top left) is equivalent to the diagram, that is obtained from diagram $L$ by replacing $L'$ with a strand without crossings.
    \end{lemma}
    \textit{Proof.} $\square$ The proof can be carried out by a sequence of equivalent diagrams shown in Fig.4. All crossings that appear in the gray area are flat. $\blacksquare$
	\subsection{Images of the map $\phi_d$}
	Let $d \geq 2$ be a natural number and let $K$ be a classical diagram on a cylinder, then there exist a natural number $m$ and diagram $K'$ with $2m$ cut ends, such that diagram $K$ is equivalent to a diagram $\fl_d(K')$ on the cylinder (Fig.5.2). Each of the $d$ rectangles on a cylinder is divided into $d$ equal parts by $d-1$ horisontal segments.
    \begin{center}
    	\includegraphics[width = 1.0\textwidth]{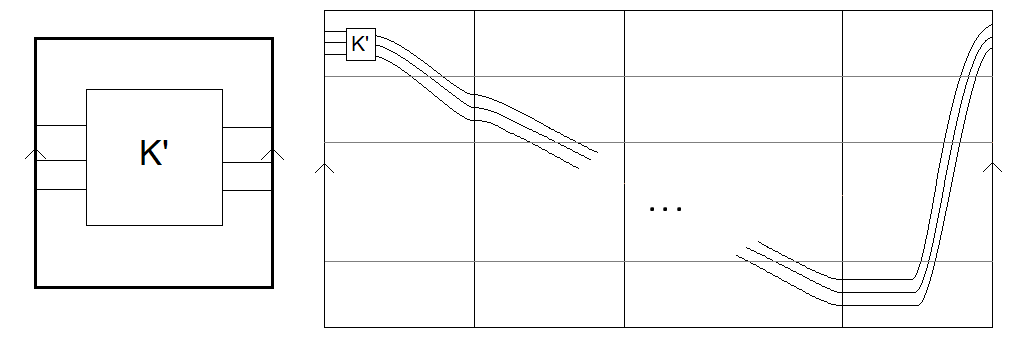}
    Figure 5.1, 5.2: Diagrams $\fl(K')$, $\fl_d(K')$ on cylinder for diagram $K'$ with $2m$ cut ends, $m = 3$
    \end{center}
    Thus $\phi_d(K)$ is equivalent to $\phi_d(\fl_d(K'))$. Fig.6.1 illustrates $\phi_d(\fl_d(K'))$ without virtual crossings. Then it can be transformed into Fig.6.2 by separating strands from each other by $2$-d and $3$-d flat moves. In Fig.6.2 each strand crosses itself $d - 1$ times.
    \begin{center}
    	\includegraphics[width = 1.0\textwidth]{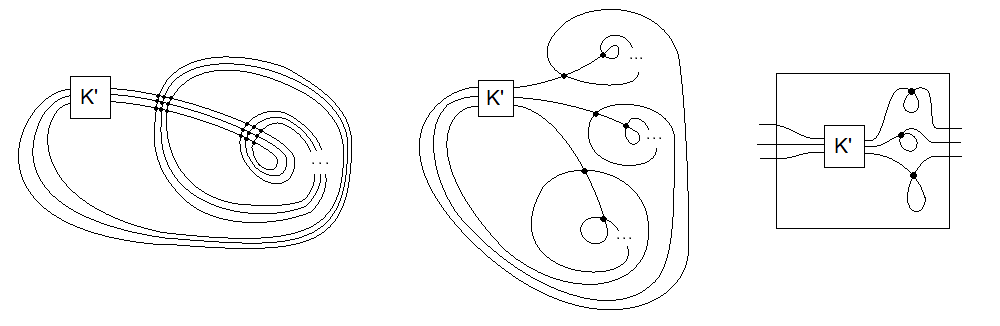}
    	Figure 6.1, 6.2: $\phi_d(\fl_d(K'))$, 6.3: $\overline{K'}$, $m = 3$
    \end{center}
    For a diagram $K'$ with $2m$ cut ends let us define diagram $\overline{K'}$ with $2m$ cut ends as showm in Fig.6.3.
    Finally Lemma 1 brings us to
    \begin{theorem} Let $K$ be a digram on a cylinder, $K'$ be a diagram with $2m$ cut ends, such that $\fl_d(K')$. Then $\phi_d(K)$ is equivalent to $C(K')$ if $d$ is odd and to $C(\overline{K'})$ else.
    \end{theorem}
    \begin{corollary}
    	For any diagram $K$ on a cylinder and natural $n \geq 2$ any smoothing of the flat-virtual diagram $\phi_d(K)$ is semi-trivial.
    \end{corollary}
    
    Calculation of $J_f \circ \phi_d$ will give sum of polynomials, multiplied by disjoint union of flat eights. For given diagram on the cylinder the calculation of $J_f \circ \phi_d$ is as follows. For every state of diagram with $2m$ cut ends reduce adjacent (from one group of $m$ ends) pairs of ends that are connected to each other by removing them in this group and connecting their opposites. Then, count number of trivial components and number of strands connecting different sides --- this will be the number of flat eight components in case when $d$ is even and number of additional trivial components otherwise. Alexander matrix of $\phi_d(K)$ (flat-virtual Alexander matrix in defined in [1]) equals to Alexander matrix of $K$, which means that for any diagram $K$ on the cylinder $\Delta(\phi_d(K)) = 0$.
    
    \textit{Example.} Calculate $J_f(\phi_4(K))$ where $K$ is a diagram on cylinder (Fig.3.2).
     For the state in Fig.7 we have two trivial components and one flat eight component which gives summand $a^{-3}a^{4 - 1}(-a^2-a^{-2})^2$\includegraphics[width=20pt]{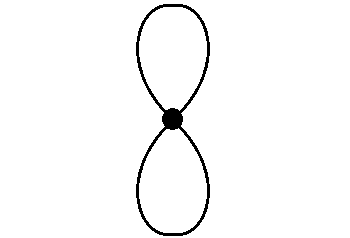} of $J_f$.
    
    \begin{center}
    	\includegraphics[width = 1.0\textwidth]{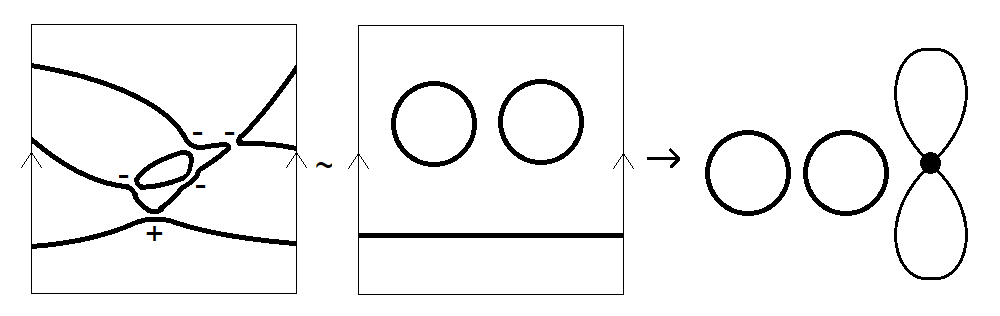}
    	Figure 7: Diagram of one summand of $J_f(\phi_4(K))$
    \end{center}
    %{\mathdfZ}
	\subsection{Images of $\phi'_{Z_{d_1} \bigoplus Z_{d_2}}$}
	Let $d_1 \geq 2$, $d_2 \geq 2$ be natural numbers. The map $\phi'_l$ makes flat virtual diagram considering angular $(\frac{2\pi k_1}{d_1}, \frac{2\pi k_2}{d_2})$ rotation, where $(k_1, k_2) \in l=Z_{d_1} \bigoplus Z_{d_2}$. For any diagram $K$ on the torus there exist natural numbers $m_1$ and $m_2$ and diagram $K'$ with $2(m_1 + m_2)$ cut ends, such that diagram $K$ is equivalent to the diagram $\fl(K')$ on the torus (Fig.9.1). Denote by $\fl_l(K')$ a special position of $\fl(K')$ that is useful for $\phi'_l$ (Fig.9.2).
	
	Thus $\phi'_l(K)$ is equivalent to $\phi'_l(\fl_l(K'))$ which contains $d_1 d_2 - 1$ blocks of $m_1 m_2$ flat, $d_2 - 1$ blocks of $m_1^2$ flat, $d_1 - 1$ blocks of $m_2^2$ flat and $d_1 d_2$ blocks of $m_1 m_2$ virtual crossings for $m_1 = 3$ and $m_2 = 2$. An example of $\phi'_l(\fl_l(K'))$ is shown in Fig.8.
	
	Denote by $D_s(m_1, m_2)$ the diagram with $2(m_1 + m_2)$ cut ends without crossings every component of which forms a "slash" and $D_b(m_1, m_2)$ a diagram with $2(m_1 + m_2)$ cut ends without crossings every component of which forms a "backslash". Notice that $\fl_l(D_s(m_1, m_2))$ is equivalent to $(m_1, m_2)$-torus link and $\fl_l(D_b(m_1, m_2))$ is equivalent to $(m_1, -m_2)$-torus link.
	
	\begin{center}
		\includegraphics[width = 0.9\textwidth]{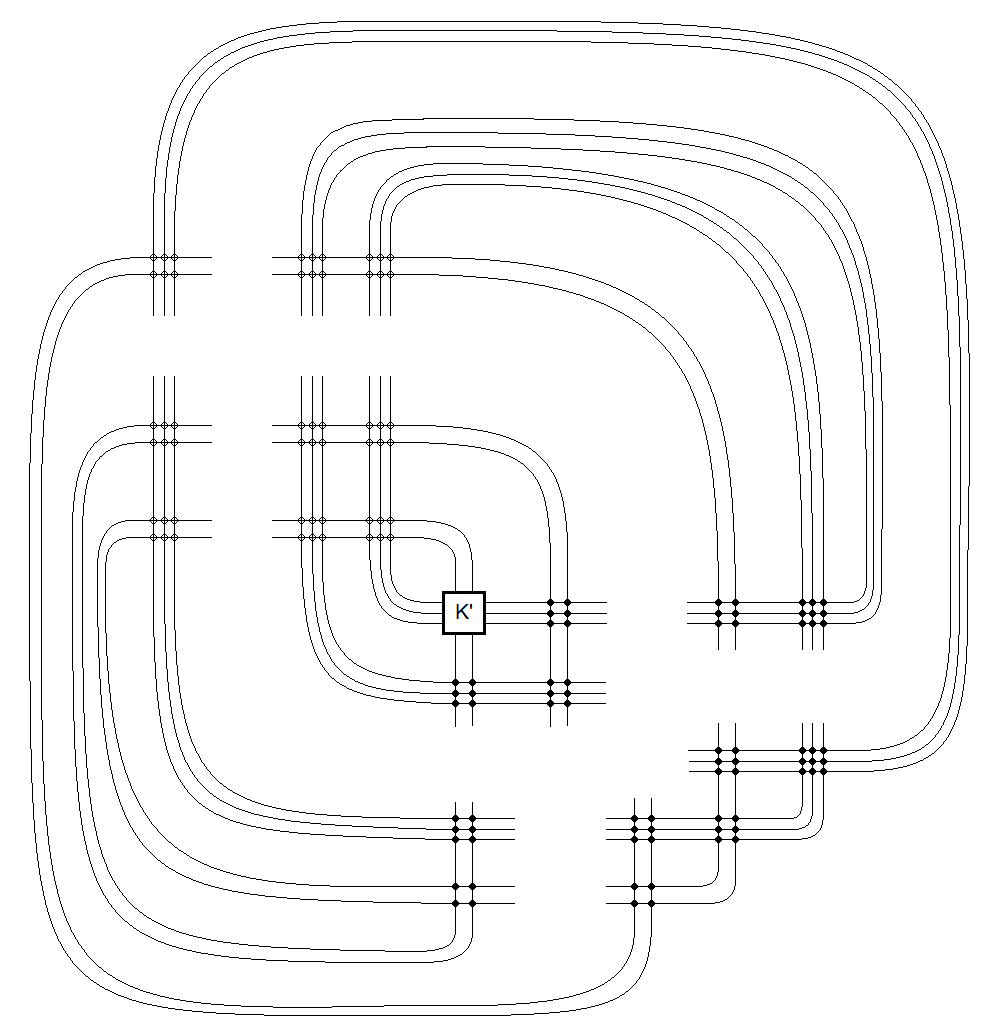}
	\end{center}
	\begin{center}
		Figure 8: $\phi'_l(\fl_l(K'))$, $m_1 = 3$, $m_2 = 2$
	\end{center}
	\begin{center}
		\includegraphics[width = 1.0\textwidth]{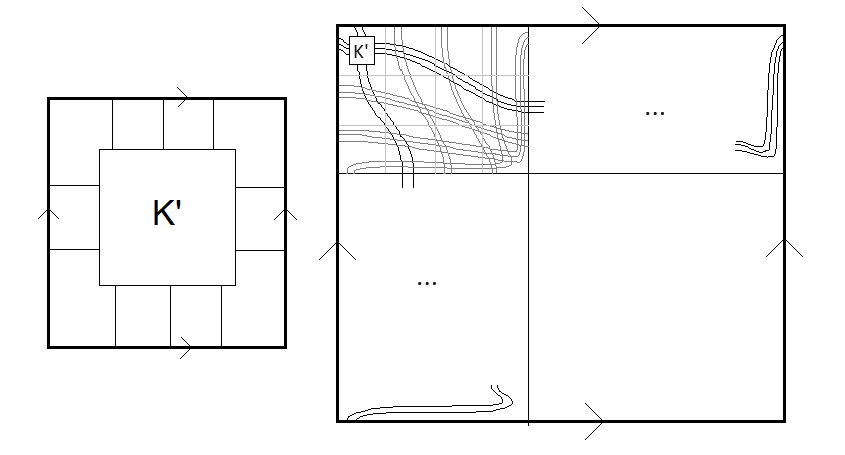}
	Figure 9.1: $\fl(K')$, $m_1 = 2$, $m_2 = 3$, 9.2: $\fl_l(K')$, $m_1=3$, $m_2=2$, $d_1=4$, $d_2=3$
	\end{center}
	\begin{center}
		\includegraphics[width = 1.0\textwidth]{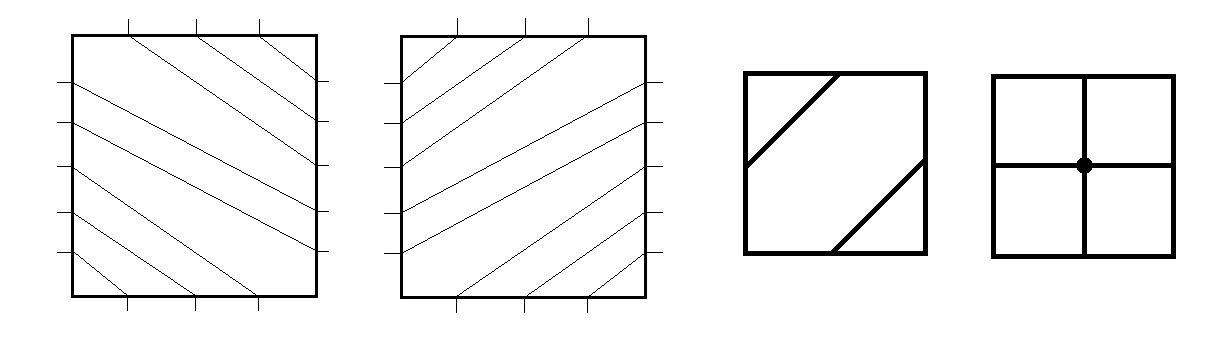}
		Figure 10.1: $D_b(5, 3)$, 10.2: $D_s(5, 3)$, 10.3, 10.4: parts of $\alpha_s(d_1, d_2, m_1, m_2)$
	\end{center}
	
	\begin{lemma}
		Fix $d_1$, $d_2$. If for any $m_1$ and $m_2$ diagrams $\phi'_l(\fl_l(D_s(m_1, m_2)))$ and $\phi'_l(\fl_l(D_b(m_1, m_2)))$ are semi-trivial, then for any $m_1$, $m_2$ and any smoothing $L'$ of diagram $K'$ with $2(m_1 + m_2)$ cut ends $\phi'_l(\fl_l(L'))$ is semi-trivial.
	\end{lemma}
	
	\textit{Proof.} $\square$ Conduct proof by induction. For $m_1 = m_2 = 1$ any smoothing of diagram $K'$ is equivalent to disjoint union of trivial components and $\fl_l(D_s(1, 1))$ or $\fl_l(D_b(1, 1))$. Now suppose lemma is true for any diagram with $\leq 2(m_1 + m_2) - 2$ ends. If $L'$ is not $D_s(m_1, m_2)$ or $D_b(m_1, m_2)$, at least one of the sides of the square border of $L'$ contains two adjacent ends, connected with each other. Then by pulling these ends (using second flat detour moves) through flat virtual part of $\fl_l(L')$ we get its equivalence to $\fl_l(\hat{L'})$ where $\hat{L'}$ is a diagram with $2(m_1 + m_2) - 2$ cut ends. $\blacksquare$
	
	\begin{theorem}
		For any $d_1$, $d_2$ and diagram $K$ on torus any smoothing of $\phi'_l(K)$ is semi-trivial.
	\end{theorem}
	
	\textit{Proof.} $\square$ Let $K'$ be a diagram with $2(m_1+m_2)$ cut ends such that $K$ is equivalent to $\fl_l(K')$. Notice that every smoothing of $\phi'_l(\fl_l(K'))$ is flat-virtually equivalent to $\phi'_l(\fl_l(L'))$, where $L'$ is a corresponding smoothing of $K'$ ($\phi'_l$ induces trivial bijection between classical crossings of $K$ and classical crossings of $\phi'_l(K)$). Thus, due to lemma 2, it is sufficient to prove that $\forall d_1, d_2, m_1, m_2$ $\phi'_l(\fl_l(D_s(m_1, m_2)))$ and $\phi'_l(\fl_l(D_b(m_1, m_2)))$ are equivalent to disjoint union of trivial and flat eight components.
	
	For numbers $d_1, d_2, m_1, m_2$ let us define two flat virtual diagrams with $2(d_1 m_2 + d_2 m_1)$ cut ends. Call them $\alpha_b(d_1, d_2, m_1, m_2)$ and $\alpha_s(d_1, d_2, m_1, m_2)$. Consider some rectangle on a plane with a grid that divides it into $d_1 m_2$ rows and $d_2 m_1$ columns of small rectangles. Let us number the rows from top to bottom starting with 0. Then the part of $\alpha_s(d_1, d_2, m_1, m_2)$ inside each small rectangle in the most right column  and in the row with number less then $\gcd(d_1 m_2, d_2 m_1)$ and relatively prime with $\frac{\gcd(d_1 m_2, d_2 m_1)}{\gcd(m_2, m_1)}$ is two segments connecting the middles of the sides of the small rectangle as shown in Fig.10.4. As for the other small rectangles $\alpha_s(d_1, d_2, m_1, m_2)$ inside them are as shown in Fig.10.3. $\alpha_b(d_1, d_2, m_1, m_2)$ can be got from $\alpha_s(d_1, d_2, m_1, m_2)$ by vertical reflection. Example of $\alpha_s(d_1, d_2, m_1, m_2)$ is shown in Fig.11.
	
	If one component of $\alpha_s(d_1, d_2, m_1, m_2)$ is oriented, then we can find a "next" component --- the first component that is met while going from selected component along $T(\alpha_s(d_1, d_2, m_1, m_2))$ oriented according to the orientation of the selected component. It is easy to check that the relation on components "the second component can be obtained from the first one by going to the "next" component a finite number of times" is an equivalence relation and components of $\alpha_s(d_1, d_2, m_1, m_2)$ are equivalent in this sense iff the corresponding parts of $\fl(\alpha_s(d_1, d_2, m_1, m_2))$ belong to one component.
	
	Denote $\alpha'_s(d_1, d_2, m_1, m_2)$ a diagram with $2(d_1 m_2 + d_2 m_1)$ cut ends that can be got from $\alpha_s(d_1, d_2, m_1, m_2)$ by replacing all \includegraphics[width=20pt]{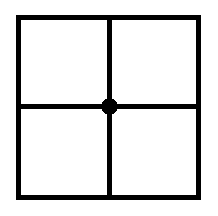} parts with \includegraphics[width=20pt]{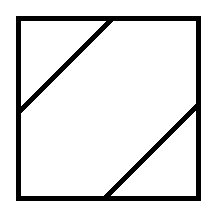}. $\alpha'_s(d_1, d_2, m_1, m_2)$ is equivalent to $D_s(d_1 m_2, d_2 m_1)$.
	
	\begin{center}
		\includegraphics[width = 1.0\textwidth]{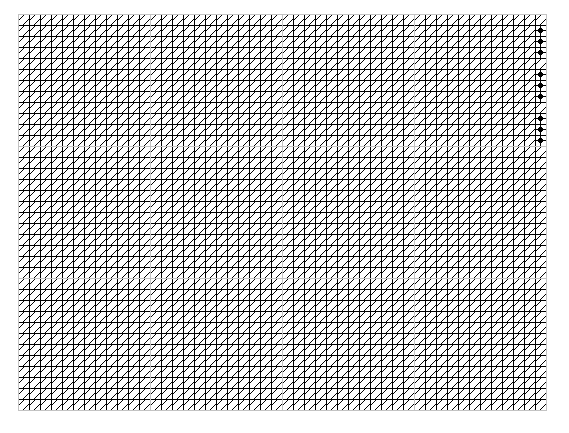}
		Figure 11: $\alpha_s(4, 8, 6, 9)$
	\end{center}
	
	$\fl(\alpha'_s(d_1, d_2, m_1, m_2))$ is $(d_1 m_2, d_2 m_1)$-torus link, which is a $\gcd(d_1 m_2, d_2 m_1)$-component link, and all its strands crossing the right side of the big rectangle in the first $\gcd(d_1 m_2, d_2 m_1)$ rows belong to pairwise distinct components. Thus, $\fl(\alpha_s(d_1, d_2, m_1, m_2))$ has $\gcd(m_1, m_2)$ components, each component crosses only itself and every group of crossings between \includegraphics[width=20pt]{uncseg.png} part in the right column belongs to the same component.
	
	Now construct curves on torus divided into $d_1 d_2$ rectangles the following way. For every component of $\fl(\alpha_s(d_1, d_2, m_1, m_2))$ choose any rectangle and an oriented component of $\alpha_s(d_1, d_2, m_1, m_2)$ and draw it in this rectangle. Then choose the "next" rectangle adjacent to the current one along the edge containing the end of drawn component. Draw in the new rectangle the component of $\alpha_s(d_1, d_2, m_1, m_2)$ that is "next" for the current component. Repeat this step until one of the two outcomes occurs:
	
	1. "Next" component has been already drawn in any of the previous rectangles.
	
	2. "Next" component is going to cross the components that have been already drawn in the "next" rectangle.
	
	In fact, case 2 never occurs and case 1 occurs only if this previous rectangle coincides with the "next" rectangle. Let us ascertain this fact. Call {\em \textbf{block}} each square block of $\gcd(d_1 m_2, d_2 m_1)^2$ small rectangles that each of the $d_1 d_2$ rectangles on torus can be divided in. The number of Blocks through which the curve passes through the left side before returning to the original {\em \textbf{block}} equals
	\begin{equation}
	    \frac{\frac{d_1 m_2 d_2}{\gcd(d_1 m_2, d_2 m_1)}\frac{d_2 m_1 d_1}{\gcd(d_1 m_2, d_2 m_1)}}{\gcd(\frac{d_1 m_2 d_1}{\gcd(d_1 m_2, d_2 m_1)}, \frac{d_2 m_1 d_2}{\gcd(d_1 m_2, d_2 m_1)})}=\frac{d_1 m_2 d_2 m_1}{\gcd(d_1 m_2, d_2 m_1) \gcd(m_1, m_2)}
    \end{equation},
    whereas $\fl(\alpha_s(d_1, d_2, m_1, m_2))$ returns to the initial Block every
    \begin{equation}
    	\frac{\frac{d_1 m_2}{\gcd(d_1 m_2, d_2 m_1)}\frac{d_2 m_1}{\gcd(d_1 m_2, d_2 m_1)}}{\gcd(\frac{d_1 m_2}{\gcd(d_1 m_2, d_2 m_1)}, \frac{d_2 m_1}{\gcd(d_1 m_2, d_2 m_1)})}=\frac{d_1 m_2 d_2 m_1}{\gcd(d_1 m_2, d_2 m_1)^2}
    \end{equation}
    {\em \textbf{blocks}}. Finally, for every {\em \textbf{block}} diagram $\fl(\alpha_s(d_1, d_2, m_1, m_2))$ goes through its left side $\frac{\gcd(d_1 m_2, d_2 m_1)}{\gcd(m_1, m_2)}$ times. Thus, every component of $\fl(\alpha_s(d_1, d_2, m_1, m_2))$ crosses the left sides of the {\em \textbf{blocks}} exactly
	\begin{equation}
	\frac{\gcd(d_1 m_2, d_2 m_1)}{\gcd(m_1, m_2)} \frac{d_1 m_2 d_2 m_1}{\gcd(d_1 m_2, d_2 m_1)^2}=\frac{d_1 m_2 d_2 m_1}{\gcd(d_1 m_2, d_2 m_1) \gcd(m_1, m_2)}
	\end{equation}
	times which means that the constructed curves do not cross themselves and all rectangles on the torus overlap on each other without coincidence of curves (only transversal intersections). Thus, these curves constitute a $(\frac{d_1 m_2}{d_1}, \frac{d_2 m_1}{d_2})$-torus link, which is equivalent to $\fl_l(D_s(m_1, m_2))$, and $\phi'_l(\fl_l(D_s(m_1, m_2)))$ is equivalent to $\T(\alpha_s(d_1, d_2, m_1, m_2))$. $\T(\alpha_s(d_1, d_2, m_1, m_2))$ consists of $\gcd(m_1, m_2)$ components, each looks like the one in Fig.12 and has $\frac{\gcd(d_1 m_2, d_2 m_1)}{\gcd{m_1, m_2}} - 1$ flat crossings. Then applying Lemma 1 we get that all components of $\phi'_l(\fl_l(D_s(m_1, m_2)))$ are either trivial or flat eight. Analogously we get that $\phi'_l(\fl_l(D_b(m_1, m_2)))$ is equivalent to $\T(\alpha_b(d_1, d_2, m_1, m_2))$ and all components of $\phi'_l(\fl_l(D_b(m_1, m_2)))$ are either trivial or flat eight. $\blacksquare$
	
	\begin{center}
		\includegraphics[width = 1.0\textwidth]{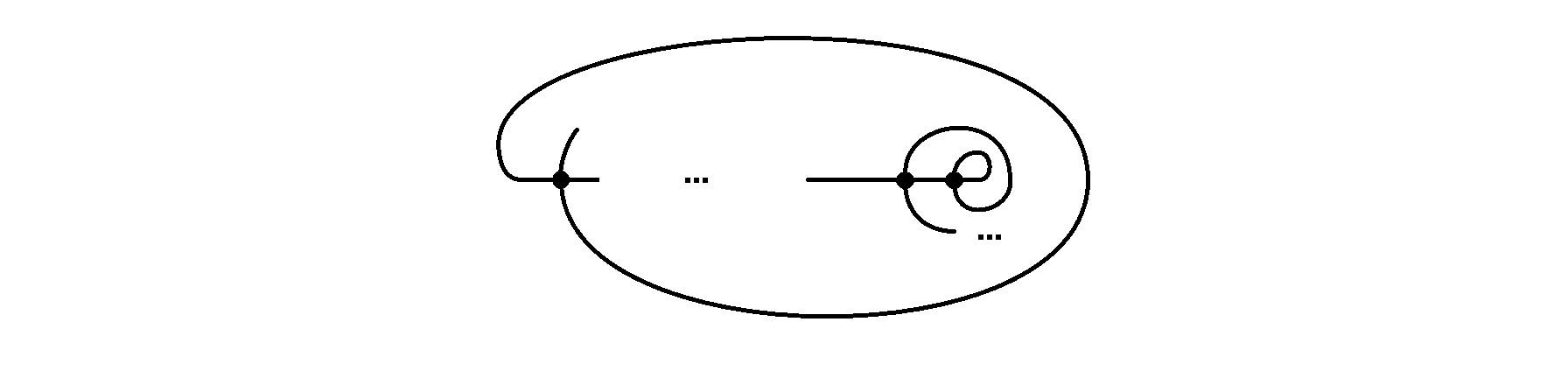}
		Figure 12: Component of $\T(\alpha_s(d_1, d_2, m_1, m_2))$
	\end{center}
	
	For every state of diagram with $2(m_1 + m_2)$ cut ends reduce adjacent (from one group of $m_1$ or $m_2$ ends) pairs of ends that are connected with each other by removing them in this group and connecting their opposites. In the end we get disjoint union of trivial components and $D_b(m'_1, m'_2)$ or $D_s(m'_1, m'_2)$ for some $m'_1 \leq m_1$ and $m'_2 \leq m_2$. Then, count number of trivial components and in case $m'_1 \neq 0$ or $m'_2 \neq 0$ calculate $\gcd(m'_1, m'_2)$ --- this will be the number of additional trivial components if $\frac{\gcd(d_1 m_2, d_2 m_1)}{\gcd{m_1, m_2}}$ is odd and the number of flat eight components else.
	
	\textit{Example.} Calculate $J_f(\phi'_l(K))$ for $d_1 = 3, d_2 = 5$ where $K$ is a diagram on torus (Fig.3.3).
	For the state in Fig.13 left we have $m'_1 = 1$, $m'_2 = 0$ and two trivial components, which means summand $a^{12} a^{3-3}(-a^2-a^{-2})^2$ of $J_f$.
	\begin{center}
		\includegraphics[width = 1.0\textwidth]{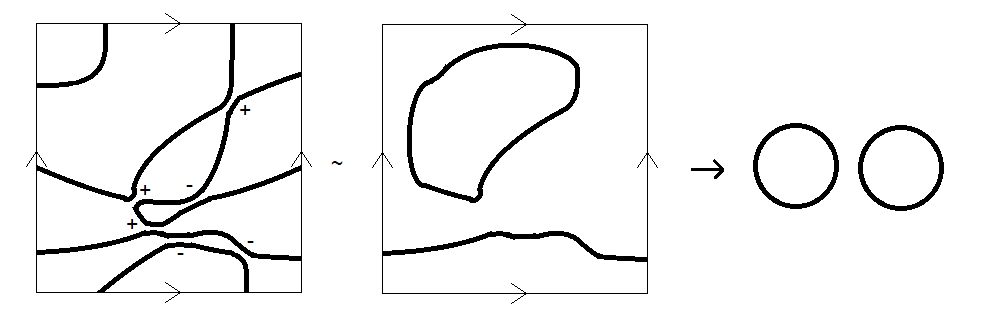}
		Figure 13: Diagram of the summand of $J_f(\phi_{Z_3 \bigoplus Z_5}(K))$
	\end{center}

\section{Acknowledgements}
    I am really grateful to V.O.Manturov and I.M.Nikonov for useful discussions.
\section{References}
	[1] L.H.Kauffman, Virtual Knot Theory, European Journal of Combinatorics, 20 (1999) 663–691.
	
	[2] V.O.Manturov, D.A.Fedoseev, S.Kim, I.M.Nikonov, Invariants and Pictures, World Scientific, 2020.
	
	[3] V.O.Manturov, D.P.Ilyutko, Virtual knots. The state of the art, World Scientific, 2012.
	
	[4] V.O.Manturov, I.M.Nikonov, Flat-virtual knot: introduction and some invariants.
	
	[5] V.O.Manturov, I.M.Nikonov, Maps from knots in the cylinder to flat-virtual knots, arxiv:2210:09689.
\end{document}